\overfullrule=0pt
\centerline {\bf A further multiplicity result for Lagrangian systems of relativistic oscillators}\par
\bigskip
\bigskip
\centerline {BIAGIO RICCERI}\par
\bigskip
\bigskip
\centerline {\it Dedicated to Professor Adrian Petru\c{s}el on his 60th birthday}
\bigskip
\bigskip
{\bf Abstract.} This is our third paper, after [4] and [5], about a joint application of the theory developed by Brezis and Mawhin in [1] with our minimax theorems ([2], [3]) to get
multiple solutions of problems of the type
$$\cases{(\phi(u'))'=\nabla_xF(t,u) & in $[0,T]$\cr & \cr
u(0)=u(T)\ , \hskip 3pt u'(0)=u'(T)\cr}$$
which are global minima of a suitable functional over a set of Lipschitzian functions. A challenging conjecture is also formulated.
 \par
\bigskip
\bigskip
{\bf 1. Introduction}\par
\bigskip
In what follows, $L, T$ are two fixed positive numbers. For each $r>0$, we set $B_r=\{x\in {\bf R}^n :|x|<r\}$ ($|\cdot|$ being the Euclidean norm on ${\bf R}^n$)
and $\overline {B_r}$ is the closure of $B_r$. The scalar product on ${\bf R}^n$ is denoted by $\langle\cdot,\cdot\rangle$.
\smallskip
We denote by ${\cal A}$ the family of all homeomorphisms $\phi$ from $B_L$ onto ${\bf R}^n$ such that $\phi(0)=0$ and $\phi=\nabla\Phi$,
where the function $\Phi:\overline {B_L}\to ]-\infty,0]$ is continuous and strictly convex in $\overline {B_L}$, and of class $C^1$ in $B_L$.
Notice that $0$ is the unique global minimum of $\Phi$ in $\overline {B_L}$.
\par
\smallskip
We denote by ${\cal B}$ the family of all functions $F:[0,T]\times {\bf R}^n\to {\bf R}$ which are measurable in $[0,T]$, of class $C^1$ in
${\bf R}^n$ and such that $\nabla_xF$ is measurable in $[0,T]$ and, 
for each $r>0$, one has $\sup_{x\in B_r}|\nabla_x F(\cdot,x)|\in L^1([0,T])$, with $F(\cdot,0)\in L^1([0,T])$. Clearly, ${\cal B}$ is
a linear subspace of ${\bf R}^{[0,T]\times {\bf R}^n}$.
\smallskip
Given $\phi\in {\cal A}$ and $F\in {\cal B}$, we consider the problem
$$\cases{(\phi(u'))'=\nabla_xF(t,u) & in $[0,T]$\cr & \cr
u(0)=u(T)\ , \hskip 3pt u'(0)=u'(T)\ .\cr}\eqno{(P_{\phi,F})}$$
A solution of this problem is any function $u:[0,T]\to {\bf R}^n$ of class $C^1$, with $u'([0,T])\subset
B_L$, $u(0)=u(T)$, $u'(0)=u'(T)$,  such that the composite function $\phi\circ u'$ is absolutely
continuous in $[0,T]$ and one has $(\phi\circ u')'(t)=\nabla_xF(t,u(t))$ for a.e. $t\in [0,T]$.\par
\smallskip
Now, we set
$$K=\{u\in \hbox {\rm Lip}([0,T],{\bf R}^n) : |u'(t)|\leq L\hskip 5pt for\hskip 3pt a.e.\hskip 3pt t\in [0,T] , u(0)=u(T)\}\ ,$$
Lip$([0,T],{\bf R}^n)$ being the space of all Lipschitzian functions from $[0,T]$ into ${\bf R}^n$.\par
\smallskip
Clearly, one has
$$\sup_{[0,T]}|u|\leq LT + \inf_{[0,T]}|u| \eqno{(1.1)}$$
for all $u\in K$. 
\smallskip
Next, consider the functional $I:K\to {\bf R}$ defined by
$$I(u)=\int_0^T(\Phi(u'(t))+F(t,u(t)))dt$$
for all $u\in K$.\par
\smallskip
In [1], Brezis and Mawhin proved the following result:\par
\medskip
THEOREM 1.1 ([1], Theorem 5.2). - {\it Any global minimum of $I$ in $K$ is a solution of problem $(P_{\phi,F})$.}\par
\medskip
On the other hand, very recently, in [6], we established the following:\par
\medskip
THEOREM 1.2 ([6], Theorem 2.2). - {\it Let $X$ be a topological space, let $E$ be a real normed space,
 let $I:X\to {\bf R}$, let $\psi:X\to E$ and let $S\subseteq E^*$ be a convex set
weakly-star dense in $E^*$.  Assume that $\psi(X)$  is not convex and that $I+\eta\circ\psi$ is lower semicontinuous and inf-compact
for all $\eta\in S$.\par
Then, there exists $\tilde\eta\in S$ such that the
function $I+\tilde\eta\circ\psi$ has at least two global minima in $X$.}\par
\medskip
The aim of this paper is to establish a new multiplicity result for the solutions of problem $(P_{\phi,F})$ as a joint application of
Theorems 1.1 and 1.2.\par
\smallskip
Notice that [4] and [5] are the only previous papers on multiple solutions for problem $(P_{\phi,F})$ which are global minima of $I$ in $K$.
\par
\bigskip
{\bf 2. The result}\par
\bigskip
Here is our result:\par
\medskip
THEOREM 2.1. - {\it Let $\phi\in {\cal A}$, $F, G\in {\cal B}$ and 
$H\in C^1({\bf R}^n)$. Assume that:\par
\noindent
$(a_1)$\hskip 5pt there exists $q>0$ such that
$$\lim_{|x|\to +\infty}{{\inf_{t\in [0,T]}F(t,x)}\over {|x|^q}}=+\infty$$
and
$$\limsup_{|x|\to +\infty}{{\sup_{t\in [0,T]}|G(t,x)|+|H(x)|}\over {|x|^q}}<+\infty\ ;$$
$(a_2)$\hskip 5pt there are $\gamma\in \{\inf_{{\bf R}^n}H,\sup_{{\bf R}^n}H\}$,
 with $H^{-1}(\gamma)$ at most countable, and
$v,w\in H^{-1}(\gamma)$ such that $\int_0^TG(t,v)dt\neq \int_0^TG(t,w)dt$.\par
Then, for each $\alpha\in L^{\infty}([0,T])$ having a constant sign and with $\hbox {\rm meas}(\alpha^{-1}(0))=0$,
there exists $(\tilde\lambda,\tilde\mu)\in {\bf R}^2$ such that the problem
$$\cases{(\phi(u'))'=\nabla_x\left (F(t,u)+\tilde\lambda G(t,u)+\tilde\mu\alpha(t)H(u)\right ) & in $[0,T]$\cr & \cr
u(0)=u(T)\ , \hskip 3pt u'(0)=u'(T)\cr}\eqno{(P)}$$
has at least two solutions which are global minima in $K$ of the functional
$$u\to \int_0^T(\Phi(u'(t))+F(t,u(t))+\tilde\lambda G(t,u(t))+\tilde\mu\alpha(t)H(u(t)))dt\ .$$}\par
\smallskip
PROOF. Fix $\alpha\in L^{\infty}([0,T])$ having a constant sign and with $\hbox {\rm meas}(\alpha^{-1}(0))=0$.
Let $C^0([0,T], {\bf R}^n)$ be the space of all continuous functions from $[0,T]$ into ${\bf R}^n$, with the norm $\sup_{[0,T]}|u|$.
We are going to apply Theorem 1.2 taking $X=K$, regarded as a subset of $C^0([0,T], {\bf R}^n)$ with the relative topology, $E={\bf R}^2$
and $I:K\to {\bf R}$, $\psi:K\to {\bf R}^2$ defined by
$$I(u)=\int_0^T(\Phi(u'(t))+F(t,u(t)))dt\ ,$$
$$\psi(u)=\left ( \int_0^TG(t,u(t))dt,\int_0^T\alpha(t)H(u(t))dt\right )$$
for all $u\in K$. Fix $(\lambda,\mu)\in {\bf R}^2$. 
 By Lemma 4.1 of [1], the function 
$I(\cdot)+\langle \psi(\cdot),(\lambda,\mu)\rangle$  is
lower semicontinuous in $K$. Let us show that it is inf-compact too. First, observe that if $P\in {\cal B}$ then, for each $r>0$, there is
$M\in L^1([0,T])$ such that
$$\sup_{x\in B_r}|P(t,x)|\leq M(t) \eqno{(2.1)}$$
for all $t\in [0,T]$. Indeed, by the mean value theorem, we have
$$P(t,x)-P(t,0)=\langle \nabla_x P(t,\xi),x\rangle$$
for some $\xi$ in the segment joining $0$ and $x$. Consequently, for all $t\in [0,T]$ and $x\in B_r$, by the Cauchy-Schwarz inequality,
we clearly have
$$|P(t,x)|\leq r\sup_{y\in B_r}|\nabla_x P(t,y)|+|P(t,0)|\ .$$
So, to get $(2.1)$,
we can choose $M(t):=r\sup_{y\in B_r}|\nabla_x P(t,y)|+|P(t,0)|$ which is in $L^1([0,T])$ since $P\in {\cal B}$. Now, by $(a_1)$,
we can fix $c_1, \delta>0$ so that
$$|G(t,x)|+|H(x)|\leq c_1|x|^q \eqno{(2.2)}$$
for all $(t,x)\in [0,T]\times ({\bf R}^n\setminus B_{\delta})$.  Then, set
$$c_2:=c_1\max\left\{|\lambda|, |\mu|\|\alpha\|_{L^{\infty}([0,T])}\right\}$$
and, by $(a_1)$ again, fix $c_3>c_2$ and $\delta_1>\delta$ so that
$$F(t,x)\geq c_3|x|^q \eqno{(2.3)}$$
for all $(t,x)\in [0,T]\times ({\bf R}^n\setminus B_{\delta_1})$. On the other hand, for what remarked above, there is
$M\in L^1([0,T])$ such that
$$\sup_{x\in B_{\delta_1}}(|F(t,x)|+|\lambda G(t,x)|+|\mu\alpha(t)H(x)|)\leq M(t) \eqno{(2.4)}$$
for all $t\in [0,T]$. Therefore, from $(2.2)$, $(2.3)$ and $(2.4)$, we infer that
$$F(t,x)\geq c_3|x|^q-M(t) \eqno{(2.5)}$$
and
$$|\lambda G(t,x)|+|\mu\alpha(t)H(x)|\leq c_2|x|^q+M(t) \eqno{(2.6)}$$
for all $(t,x)\in [0,T]\times {\bf R}^n$. Set
$$b:=T\Phi(0)-2\int_0^TM(t)dt\ .$$
For each $u\in K$, with $\sup_{[0,T]}|u|\geq LT$, taking $(1.1)$, $(2.5)$ and $(2.6)$
into account, we have
$$I(u)+\langle \psi(u), (\lambda,\mu)\rangle\geq T\Phi(0)+\int_0^TF(t,u(t))dt-\int_0^T|\lambda G(t,u(t))|dt-
\int_0^T|\mu\alpha(t)H(u(t))|dt$$
$$\geq T\Phi(0)+c_3\int_0^T|u(t)|^qdt-\int_0^TM(t)dt-c_2\int_0^T|u(t)|^qdt-\int_0^TM(t)dt$$
$$\geq (c_3-c_2)T\inf_{[0,T]}|u|^q-b\geq (c_3-c_2)T\left (\sup_{[0,T]}|u|-LT\right )^q+b\ .$$
Consequently
$$\sup_{[0,T]}|u|\leq \left ( {{I(u)+\langle \psi(u), (\lambda,\mu)\rangle - b}\over {(c_3-c_2)T}}\right )^{1\over q}+LT\ .\eqno{(2.7)}$$
Fix $\rho\in {\bf R}$. By $(2.7)$, the set 
$$C_{\rho}:=\{u\in K : I(u)+\langle \psi(u), (\lambda,\mu)\rangle\leq \rho\}$$
 turns out to be bounded. 
Moreover, the functions belonging to $C_{\rho}$ are
equi-continuous since they lie in $K$. As a consequence, by the Ascoli-Arzel\`a theorem, $C_{\rho}$ is relatively compact in $C^0([0,T],{\bf R}^n)$.  By lower semicontinuity, $C_{\rho}$ is closed in $K$. But $K$ is closed in $C^0([0,T], {\bf R}^n)$ and hence $C_{\rho}$ is compact. The inf-compactness of $I(\cdot)+\langle \psi(\cdot), (\lambda,\mu)\rangle$ is so shown. Now, we are going to prove that the
set $\psi(K)$ is not convex. By $(a_2)$, the set $\left\{\int_0^TG(t,x)dt : x\in H^{-1}(\gamma)\right\}$ is at most countable since 
$H^{-1}(\gamma)$ is so. Hence, since $\int_0^TG(t,v)dt\neq \int_0^TG(t,w)dt$, we can fix $\lambda\in ]0,1[$ so that
$$\int_0^TG(t,x)dt\neq \int_0^TG(t,w)dt+\lambda\left (\int_0^TG(t,v)dt-\int_0^TG(t,w)dt\right ) \eqno{(2.8)}$$
for all $x\in H^{-1}(\gamma)$. Since $K$ contains the constant functions, the points
$$\left (\int_0^TG(t,v)dt,\gamma\int_0^T\alpha(t)dt\right )$$ and 
$$\left (\int_0^TG(t,w)dt,\gamma\int_0^T\alpha(t)dt\right )$$ belong to $\psi(K)$.
So, to show that $\psi(K)$ is not convex, it is enough to check that
the point $$\left (\int_0^TG(t,w)dt+\lambda\left (\int_0^TG(t,v)dt-\int_0^TG(t,w)dt\right ),\gamma\int_0^T\alpha(t)dt\right )$$
 does not belong to $\psi(K)$.
Arguing by contradiction, suppose that there exists $u\in K$ such that
$$\int_0^TG(t,u(t))dt= \int_0^TG(t,w)dt+\lambda\left (\int_0^TG(t,v)dt-\int_0^TG(t,w)dt\right ) \ ,\eqno{(2.9)}$$
$$\int_0^T\alpha(t)H(u(t))dt=\gamma\int_0^T\alpha(t)dt\ .\eqno{(2.10)}$$
Since the functions $\alpha$ and $H\circ u-\gamma$ do not change sign, $(2.10)$ implies that $\alpha(t)(H(u(t))-\gamma)=0$
a.e. in $[0,T]$. Consequently, since $\hbox {\rm meas}(\alpha^{-1}(0))=0$, we have $H(u(t))=\gamma$ a.e. in $[0,T]$ and hence
$H(u(t))=\gamma$ for all $t\in [0,T]$ since $H\circ u$ is continuous. In other words, the connected set $u([0,T])$ is contained
in $H^{-1}(\gamma)$ which is at most countable. This implies that the function $u$ must be constant and so $(2.9)$ contradicts $(2.8)$.
Therefore, $I$ and $\psi$ satisfy the assumptions of Theorem 1.2 and hence there exists $(\tilde\lambda,\tilde\mu)\in
{\bf R}^2$ such that the function $I(\cdot)+\langle \psi(\cdot), (\tilde\lambda,\tilde\mu)\rangle$ has at least two global minima
in $K$. Thanks to Theorem 1.1, they are solutions of Problem (P), and the proof is complete.\hfill $\bigtriangleup$\par
\medskip
REMARK 1. Of course, $(a_2)$ is the leading assumption of Theorem 2.1. The request that $H^{-1}(\gamma)$ must be at most
countable cannot be removed. Indeed, if we remove such a request, we could take $H=0$, $G(t,x)=\langle x,\omega\rangle$,
with $\omega\in {\bf R}^n\setminus \{0\}$ and $F(t,x)={{1}\over {p}}|x|^p$, with $p>1$. Now, observe that, by Proposition 3.2 of [1],
for all $\lambda\in {\bf R}$, the problem
$$\cases{(\phi(u'))'=|u|^{p-2}u+\lambda\omega & in $[0,T]$\cr & \cr
u(0)=u(T)\ , \hskip 3pt u'(0)=u'(T)\cr}$$
has a unique solution. To the contrary, the question of whether the condition $\int_0^TG(t,v)dt\neq \int_0^TG(t,w)dt$ (keeping $v\neq w$) can be dropped remains open at present. We feel, however, that it cannot be removed. In this connection, we propose the following\par
\medskip
CONJECTURE 2.1. - {\it There exist $\phi\in {\cal A}$, $F\in {\cal B}$,
$H\in C^1({\bf R}^n)$, $\alpha\in L^{\infty}([0,T])$, with $\alpha\geq 0$ and $\hbox {\rm meas}(\alpha^{-1}(0))=0$,
 and $q>0$ for which the following assertions hold:\par
\noindent
$(b_1)$\hskip 5pt 
$$\lim_{|x|\to +\infty}{{\inf_{t\in [0,T]}F(t,x)}\over {|x|^q}}=+\infty$$
and
$$\limsup_{|x|\to +\infty}{{|H(x)|}\over {|x|^q}}<+\infty\ ;$$
$(b_2)$\hskip 5pt the function $H$ has exactly two global minima\ ;\par
\noindent
$(b_3)$\hskip 5pt  for each $\mu\in {\bf R}$, the functional
$$u\to \int_0^T(\Phi(u'(t))+F(t,u(t))+\mu\alpha(t)H(u(t)))dt$$
has a unique global minimum in $K$.}\par
\bigskip
\bigskip
{\bf Acknowledgement.} The author has been supported by the Gruppo Nazionale per l'Analisi Matematica, la Probabilit\`a e 
le loro Applicazioni (GNAMPA) of the Istituto Nazionale di Alta Matematica (INdAM) and by the Universit\`a degli Studi di Catania, PIACERI 2020-2022, Linea di intervento 2, Progetto ”MAFANE”.
\vfill\eject
\centerline {\bf References}\par
\bigskip
\bigskip
\noindent
[1]\hskip 5pt H. BREZIS and J. MAWHIN, {\it Periodic solutions of Lagrangian systems of relativistic oscillators},
Commun. Appl. Anal., {\bf 15} (2011), 235-250.\par
\smallskip
\noindent
[2]\hskip 5pt B. RICCERI, {\it Well-posedness of constrained minimization
problems via saddle-points}, J. Global Optim., {\bf 40} (2008),
389-397.\par
\smallskip
\noindent
[3]\hskip 5pt B. RICCERI, {\it On a minimax theorem: an improvement, a new proof and an overview of its applications},
Minimax Theory Appl., {\bf 2} (2017), 99-152.\par
\smallskip
\noindent
[4]\hskip 5pt B. RICCERI, {\it Multiple periodic solutions of Lagrangian systems of relativistic oscillators}, in  
``Current Research in Nonlinear Analysis - In Honor of Haim Brezis and Louis Nirenberg", Th. M. Rassias ed., 249-258, Springer, 2018.
\par
\smallskip
\noindent
[5]\hskip 5pt B. RICCERI, {\it Another multiplicity result for the periodic solutions of certain systems}, Linear Nonlinear Anal., {\bf 5} (2019), 371-378.\par
\smallskip
\noindent
[6]\hskip 5pt B. RICCERI, {\it Multiplicity theorems involving functions with non-convex range}, submitted to Stud. Univ. Babe\c{s}-Bolyai Math.
\bigskip
\bigskip
Department of Mathematics and Informatics\par
University of Catania\par
Viale A. Doria 6\par
95125 Catania, Italy\par
{\it e-mail address:} ricceri@dmi.unict.it

\bye